\documentclass[aop,MSNbibl,seceqn,nameyear,dvips]{arximspdf}
\usepackage{graphicx}
\doi{10.1214/11-AOP652}
\volume{39}
\issue{5}
\pubyear{2011}
\firstpage{1844}
\lastpage{1863}

\makeatletter

\newcommand{\iiint}{\int\!\!\int\!\!\int}

\newcommand{\Mcardy}{\texttt{cardy} }
\newcommand{\Mcr}{\texttt{cr} }
\newcommand{\ML}{\texttt{L} }
\newcommand{\MLd}{\texttt{L1} }
\newcommand{\Mh}{\texttt{h} }

\newcommand{\cardy}{\operatorname{cardy}}
\newcommand{\watts}{\operatorname{watts}}
\newcommand{\tripod}{\operatorname{tripod}}
\newcommand{\SLEtripod}{\operatorname{SLEtripod}}
\newcommand{\cratio}{\operatorname{cr}}
\newcommand{\drift}{\operatorname{drift}}

\newtheorem{theorem}{Theorem}[section]
\newtheorem{corollary}[theorem]{Corollary}
\newtheorem{lemma}[theorem]{Lemma}

\newcommand{\R}{\mathbb{R}}
\newcommand{\C}{\mathbb{C}}
\newcommand{\N}{\mathbb{N}}

\newcommand{\eps}{\varepsilon}
\newcommand{\p}{{\partial}}

\makeatother

\begin{document}
\begin{frontmatter}

\title{Schramm's proof of Watts' formula}
\runtitle{Schramm's proof of Watts' formula}

\begin{aug}
\author[A]{\fnms{Scott} \snm{Sheffield}\ead[label=u1,url]{http://math.mit.edu/\textasciitilde sheffield}}
and
\author[B]{\fnms{David B.} \snm{Wilson}\corref{}\ead[label=e2]{dbwilson@microsoft.com}\ead[label=u2,url]{http://dbwilson.com}}
\runauthor{S. Sheffield and D. B. Wilson}
\affiliation{Massachussetts Institute of Technology and Microsoft Research}
\address[A]{Department of Mathematics\\
Massachussetts Institute of Technology\\
77 Massachusetts Ave\\
Cambridge, Massachusetts 02139\\
USA\\
\printead{u1}} 
\address[B]{Microsoft Research\\
One Microsoft Way\\
Redmond, Washington 98052\\
USA\\
\printead{u2}}
\end{aug}

\received{\smonth{3} \syear{2010}}
\revised{\smonth{2} \syear{2011}}

%
\begin{abstract}
G\'{e}rard Watts predicted a formula for the probability in
percolation that there is both a left--right and an up--down crossing,
which was later proved by Julien Dub\'{e}dat. Here we present a
simpler proof due to Oded Schramm, which builds on Cardy's formula
in a conceptually appealing way: the triple derivative of Cardy's
formula is the sum of two multi-arm densities. The relative sizes
of the two terms are computed with Girsanov conditioning. The
triple integral of one of the terms is equivalent to Watts' formula.
For the relevant calculations, we present and annotate Schramm's
original (and remarkably elegant) Mathematica code.
\end{abstract}

%
\begin{keyword}[class=AMS]
\kwd[Primary ]{60J67}
\kwd[; secondary ]{82B43}.
\end{keyword}
\begin{keyword}
\kwd{Schramm--Loewner evolution (SLE)}
\kwd{percolation}.
\end{keyword}

\end{frontmatter}

\section{Watts' formula}

\begin{figure}

\includegraphics{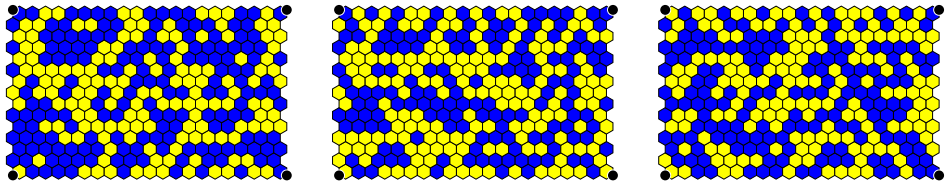}

\caption{In the left panel, there is no left--right
crossing in blue hexagons. In the second panel there is a blue
left--right crossing, but no blue up--down crossing. In the third panel,
there are both blue left--right and blue up--down crossings. Cardy's
formula gives the probability of a left--right crossing in a domain,
while Watts' formula gives the probability that there is both a
left--right crossing and an up--down crossing.} \label{fig:cardy-watts}
\end{figure}

When Langlands, Pichet, Pouliot and Saint-Aubin (\citeyear{LPPS}) were
doing computer simulations to test the conformal invariance of
percolation, there were several different events whose probability
they measured. The first event that they studied was the probability
that there is a percolation crossing connecting two disjoint boundary
segments. Using conformal field theory, \citet{Cardy}
derived his now-famous formula for this crossing probability, and the
formula was later proved rigorously by \citet{Smirnov} for site
percolation on the hexagonal lattice. The next event
that Langlands et al. tested was the probability that there
is both a percolation crossing connecting the two boundary segments
and a percolation crossing connecting the complementary boundary
segments (see Figure~\ref{fig:cardy-watts}). This probability also appeared to be conformally invariant,
but finding a~formula for it was harder, and it was not until several
years after Cardy's formula that  \citet{Watts} proposed his formula
for the probability of this double crossing. Watts
considered the derivation of the formula unsatisfactory, even by the
standards of
physics, but it matched the data of Langlands et
al. very well, which lent credibility to the formula. Watts'
formula was proved rigorously by \citet{DubedatWatts}.

%
%

To express Cardy's formula and Watts' formula for the two types of
crossing events, since the scaling limit of percolation is conformally
invariant, it is enough to give these probabilities for one canonical
domain, and this is usually taken to be the upper half-plane.
There are four points on the boundary of the domain (the real line).
Label them in increasing order $x_1$, $x_2$, $x_3$ and $x_4$. Cardy's
formula is then the
probability that there is a percolation crossing from the interval
$[x_1,x_2]$ to the interval $[x_3,x_4]$. Again by conformal invariance,
we may map the upper half-plane to itself so that $x_1\to0$, $x_3 \to1$
and $x_4\to\infty$. The remaining point $x_2$ gets mapped to
%
%
\begin{equation} s =\cratio(x_1,x_2,x_3,x_4) := \frac
{(x_2-x_1)(x_4-x_3)}{(x_3-x_1)(x_4-x_2)},
\end{equation}
which is a point in $(0,1)$ known as the cross-ratio. Both Cardy's
formula and
Watts' formula are expressed in terms of the cross-ratio.
Cardy's formula for the probability of a percolation crossing is
%
%
\begin{equation}\label{e.cardy} \cardy(s) := \frac{\Gamma(
2/3)}{\Gamma(4/3)\Gamma(1/3)}  s^{1/3}\,\hspace*{0.4pt}{}_2\hspace*{-0.4pt}F_1\biggl(\frac
13,\frac23;\frac43;s\biggr),
\end{equation}
where $\Gamma$ is the gamma function, and ${}_2\hspace*{-0.4pt}F_1$ is the hypergeometric
function defined by
\[
{}_2\hspace*{-0.4pt}F_1(a,b;c;z) = \sum_{n=0}^\infty\frac{(a)_n (b)_n}{(c)_n n!} z^n,
\]
where $a,b,c\in\C$ are parameters, $c\notin-\N$
(where $\N=\{0,1,2,\ldots\}$), and $(\ell)_n$ denotes $\ell
(\ell+1)\cdots(\ell+n-1)$. This series converges for $z\in\C$ when $|z|<1$,
and the hypergeometric function is defined by analytic continuation elsewhere
(though it is then not always single-valued).

\begin{figure}

\includegraphics{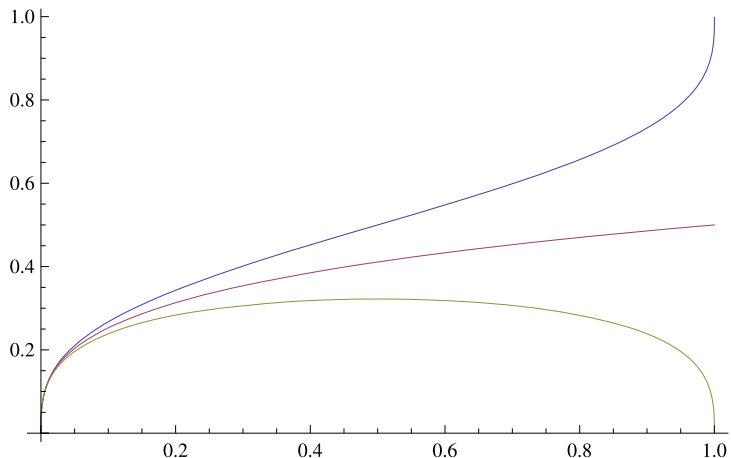}

\caption{Cardy's formula (upper curve), Watts'
formula (lower curve), and a tripod probability (defined in
Section \protect\ref{s.tripodreduction}) as a function of the
cross-ratio $s$.}
\label{fig:plot}
\end{figure}

By comparison, Watts' formula for the probability of the two crossings is
the same as Cardy's formula minus another term
%
%
\begin{eqnarray}
\label{e.watts}
\watts(s)&:=& \frac{\Gamma(2/3)}{\Gamma(4/3)\Gamma(
1/3)}   s^{1/3}\,\hspace*{0.4pt}{{}_2\hspace*{-0.4pt}F_1\biggl(\frac13,\frac23;\frac43;s\biggr)}
\nonumber\\[-8pt]\\[-8pt]
&&{}- \frac{1}{\Gamma(1/3)\Gamma(2/3)}
s\,{\hspace*{0.4pt}{}_3\hspace*{-0.4pt}F_2\biggl(1,1,\frac43;2,\frac53;s\biggr)},\nonumber
\end{eqnarray}
where ${}_3\hspace*{-0.4pt}F_2$ is the generalized hypergeometric function.
The functions $\operatorname{cardy}(s)$ and $\operatorname{watts}(s)$
are shown in Figure \ref{fig:plot}.
[The reader should not be intimidated by these formulae; the parts of
the proof involving hypergeometric functions can be handled mechanically
with the aid of Mathematica.
See also \citet{Watts} and \citet{Maier}
for equivalent double-integral formulations of Watts' formula.]

%
%

Schramm thought that Dub\'{e}dat's paper on Watts' formula was an
exciting development and started reading it as soon as it appeared in
the arXiv. Schramm sometimes presented papers to interested people at
Microsoft Research: for example, he presented Smirnov's proof of
Cardy's formula when it came out [\citet{Smirnov}],
as well as Dub\'{e}dat's paper on Watts' formula [\citet{DubedatWatts}],
and later Zhan's paper on the
reversibility of SLE$_\kappa$ for $\kappa\leq4$ [\citet{Zhan}]. In the
course of
reaching his own understanding of Watts' formula, Schramm simplified
Dub\'{e}dat's proof, with the help of a Mathematica notebook, and it was
this version that he presented at Microsoft on May 17, 2004. This
proof did not come up again until an August 2008
Centre de Recherches Math\'{e}matiques (CRM) meeting on SLE in
Montr\'{e}al, after a talk by Jacob Simmons on his work with Kleban and
Ziff on ``Watts' formula and logarithmic conformal field theory''
[\citet{SKZ}]. Schramm mentioned that he had an easier proof of Watts'
formula, which he recalled after just a few minutes. The people who
saw his version of the proof thought it was very elegant and strongly
encouraged him to write it up. The next day Oded wrote down an
outline of the proof, but he tragically died a few weeks later. There
is interest in seeing a written version of Schramm's version of the
proof, so we present it here.

\section{Outline of proof}
This is a slightly edited version of the proof outline that Oded wrote
down at the CRM. Steps 1 and 2 are the same as in Dub\'{e}dat's proof,
but with step 3 the proofs diverge. We will expand on these steps of
the outline (with slightly modified notation) in subsequent sections.

\begin{itemize}
\item Reduce to the problem of calculating the probability that there
is a crossing up--down
which also connects to the right.
\item Further reduce to the following problem. In the upper half-plane,
say, mark points
$-\infty<y_1<x_0<y_2<y_3=\infty$. Let $\gamma$ be the SLE$_6$
interface started from~$x_0$.
Let $\tau:=\inf\{t\ge0\dvtx\gamma_t\in\R\setminus[y_1,y_2]$ and
$\sigma:=\sup\{t<\tau\dvtx\gamma_t\in\R\}$.
Calculate $\mathbf{P}[\gamma_\tau\in[y_2,y_3], \gamma_\sigma\in
[x_0,y_2]]$.
\item Let $\sigma_1:=\sup\{t<\tau\dvtx\gamma_t\in[y_1,x_0]\}$
and $\sigma_2:=\sup\{t<\tau\dvtx\gamma_t\in[x_0,y_2]\}$.
Now calculate the probability density of the event
$\gamma_{\sigma_1}=z_1,\gamma_{\sigma_2}=z_2,\gamma_\tau=z_3$
as $h(z_1,x_0,z_2,z_3):=\p_{z_1}\p_{z_2}\p_{z_3}\operatorname
{Cardy}(z_1,x_0,z_2,z_3)$.
\item Now, $h$ [times certain derivatives] is a martingale for the
corresponding diffusion. Consider the
Doob-transform ($h$-transform) of the diffusion with this $h$. This
corresponds to
conditioning on this probability zero event. For the Doob-transform, calculate
the probability that $\sigma_2>\sigma_1$. This comes out to be a
hypergeometric
function $g$. Finally,
\[
\operatorname{Watts}(y_1,x_0,y_2,y_3)=
\int_{[y_1,x_0]}dz_1\int_{[x_0,y_2]} dz_2 \int_{[y_2,y_3]} g h\,
dz_3,
\]
(or more precisely, the three-arm probability),
and use integration by parts.
\end{itemize}

\section{Reduction to tripod probabilities} \label{s.tripodreduction}

$\!\!\!$The initial reduction, which is step~1 of the proof, has been derived
by multiple people independently. The first place that it appeared in
print appears to be in Dub\'{e}dat's (\citeyear{dubedatBM}) paper, where it
is credited to
Werner, who, in turn, is sure that it must have been known earlier.
In the interest of keeping the exposition self-contained, we explain
this reduction.

It is an elementary fact that exactly one of the following two events occurs:
\begin{longlist}[(2)]
\item[(1)]
there is a horizontal blue crossing in the rectangle (i.e., a path of
blue hexagons connecting the left
and right edges of the rectangle), which we denote by $H_b$;
\item[(2)] there is a vertical yellow crossing, which we denote by $V_y$.
\end{longlist}

\begin{center}

\includegraphics{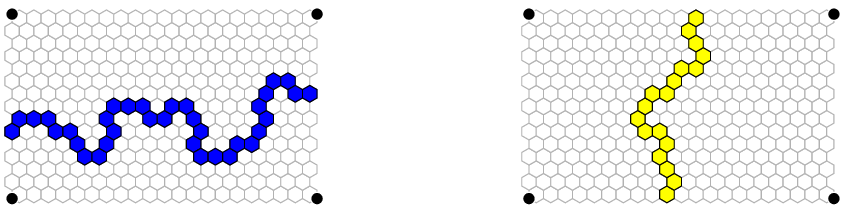}
\vspace*{3pt}
\end{center}

If there is a horizontal crossing, then by considering the region
beneath it, using the above fact, either it connects to the bottom
edge (forming a T shape) or else there is another crossing beneath it
of the opposite color. Since there are finitely many hexagons, there
must be a bottom-most crossing, which then necessarily forms a T shape.
Thus exactly one of the following three events occurs:
\begin{longlist}[(2)]
\item[(1)] there is no horizontal crossing of either color (denoted by $N$);
\item[(2)] there is a blue T (denoted $T_b$);
\item[(3)] there is a yellow T (denoted $T_y$).
\end{longlist}

\begin{center}

\includegraphics{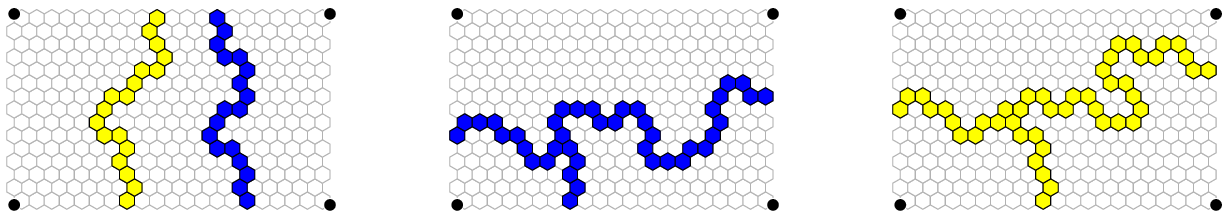}
\vspace*{3pt}
\end{center}

Of course the latter two events have equal probability, so we have
\[
\Pr[N] + 2 \Pr[T_b] = 1.
\]

Recall again that there is either a blue horizontal crossing or a
yellow vertical crossing but not both. We can decompose the yellow
vertical crossing event into two subevents according to whether or not
there is also a~yellow horizontal crossing. The first subevent is,
of course, the event we are interested in (with blue and yellow
reversed), and the second subevent is identical to the event $N$.\vspace*{3pt}

\begin{center}

\includegraphics{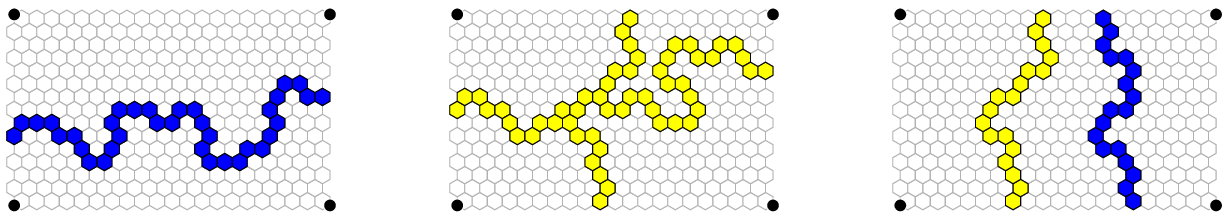}
\end{center}

Thus we have
\[
\Pr[H_b] + \Pr[H_y \wedge V_y] + \Pr[N] = 1.
\]
Combining these equations, we see that
\[
\Pr[H_b \wedge V_b] = 2 \Pr[T_b] - \Pr[H_b].
\]
In the limit of large grids with cross ratio $s$, the third term is
given by Cardy's formula, $\cardy(s)$,
and we seek to show that the left-hand side is given by Watts'
formula, $\watts(s)$. Let us give another name for
what we expect to be the limit of the second term. Define $\tripod(s)$
to satisfy
\[
\watts(s)=2\tripod(s) - \cardy(s),
\]
that is [substituting (\ref{e.cardy}) and (\ref{e.watts})],
\begin{eqnarray*}\tripod(s) &=& \frac{\watts(s) + \cardy(s)}{2} \\
&=& \frac{\Gamma(2/3)}{\Gamma(4/3)\Gamma(1/3)}
s^{1/3}\,\hspace*{0.4pt}{{}_2\hspace*{-0.4pt}F_1\biggl(\frac13,\frac23;\frac43;s\biggr)}
\\
&&{} - \frac{1}{2\Gamma(1/3)\Gamma(2/3)}
s\,\hspace*{0.4pt}{{}_3\hspace*{-0.4pt}F_2\biggl(1,1,\frac43;2,\frac53;s\biggr)}.
\end{eqnarray*}
Then in light of Cardy's formula, proving Watts' formula is equivalent
to showing that $\Pr[T_b]$ is given by $\tripod(s)$ in the fine mesh limit.\vadjust{\goodbreak}

\section{Comparison with SLE$_6$}
\subsection{Discrete derivatives of the tripod probability}

\begin{figure}

\includegraphics{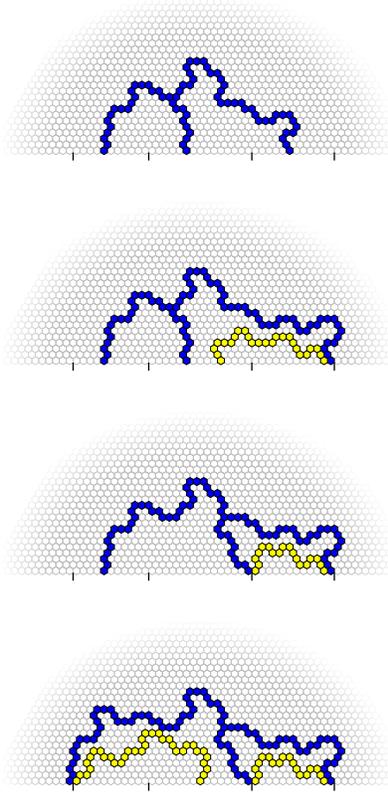}

\caption{The discrete triple partial derivative of the
tripod probability is the probability of a~multi-arm event. The top
panel illustrates the event whose probability is $P_T[x_1,x_2,x_3,x_4]$,
the next panel illustrates $\Delta_{x_4} P_T[x_1,x_2,x_3,x_4]$, the
third panel illustrates $\Delta_{x_3} \Delta_{x_4}
P_T[x_1,x_2,\allowbreak x_3,x_4]$ and the bottom panel illustrates
$-\Delta_{x_1} \Delta_{x_3} \Delta_{x_4} P_T[x_1,x_2,x_3,x_4]$.}
\label{fig:dddtripod}
\end{figure}

Consider percolation on the upper half-plane triangular lattice, and
let $P_T[x_1,x_2,x_3,x_4]$ be the probability of a blue tripod
connecting the intervals $(x_1,x_2)$ and $(x_2,x_3)$ and $(x_3,x_4)$
when the four (here discrete) locations are $x_1 < x_2<x_3<x_4$ (each
of which is a point between two boundary hexagons; see the upper image
in Figure \ref{fig:dddtripod}).

%
%

Then $\Delta_{x_4} P_T[x_1,x_2,x_3,x_4] := P_T[x_1,x_2,x_3,x_4] -
P_T[x_1,x_2,x_3,x_4-1]$ gives the probability that there is a crossing
tripod for $(x_1,x_2,x_3,x_4)$ but \textit{not} for
$(x_1,x_2,x_3,x_4-1)$. (Here we assume that the lattice spacing is $1$.)
Since the crossing tripod for
$(x_1,x_2,x_3,x_4)$ does not extend to a crossing tripod for
$(x_1,x_2,x_3,x_4-1)$, there must be a path of the opposite color from
the hexagon just to the left of $x_4-1$ to the interval between
$x_2+1$ and $x_3+1$; this event is represented by the second image in
Figure \ref{fig:dddtripod}. Similarly,
\[
-\Delta_{x_1} \Delta_{x_3} \Delta_{x_4}P_T[x_1,x_2,x_3,x_4]
\]
gives
the probability of a multi-arm event such as the one in the bottom
image in Figure \ref{fig:dddtripod}.

By summing these discrete differences, it is straightforward to write
\[
P_T[x_1,x_2,x_3,x_4] = \sum_{c\in(x_3,x_4]} \sum_{b\in(x_2,x_3]}
\sum_{a\in(x_1,x_2]}
-\Delta_{x_1} \Delta_{x_3} \Delta_{x_4}P_T[a,x_2,b,c].
\]
If there is a blue tripod connecting the intervals $(x_1,x_2)$,
$(x_2,x_3)$ and $(x_3,x_4)$, then there is only one cluster containing
such a tripod.
This formula can be interpreted as partitioning the tripod event into
multi-arm events of the type shown in the bottom panel
of Figure \ref{fig:dddtripod}.
The triple $(a,b,c) \in(x_1,x_2] \times(x_2,x_3] \times(x_3,x_4]$ is
uniquely determined by the tripod:
$a$ is (half a lattice spacing to the right of) the rightmost boundary
point of the tripod cluster in the interval $(x_1,x_2)$, $b$ is (just
right of) the rightmost point of the tripod cluster in $(x_2,x_3)$ and
$c$ is (just right of) the leftmost point of the tripod cluster in $(x_3,x_4)$.

\subsection{Discrete derivatives of the crossing probability}

\begin{figure}

\includegraphics{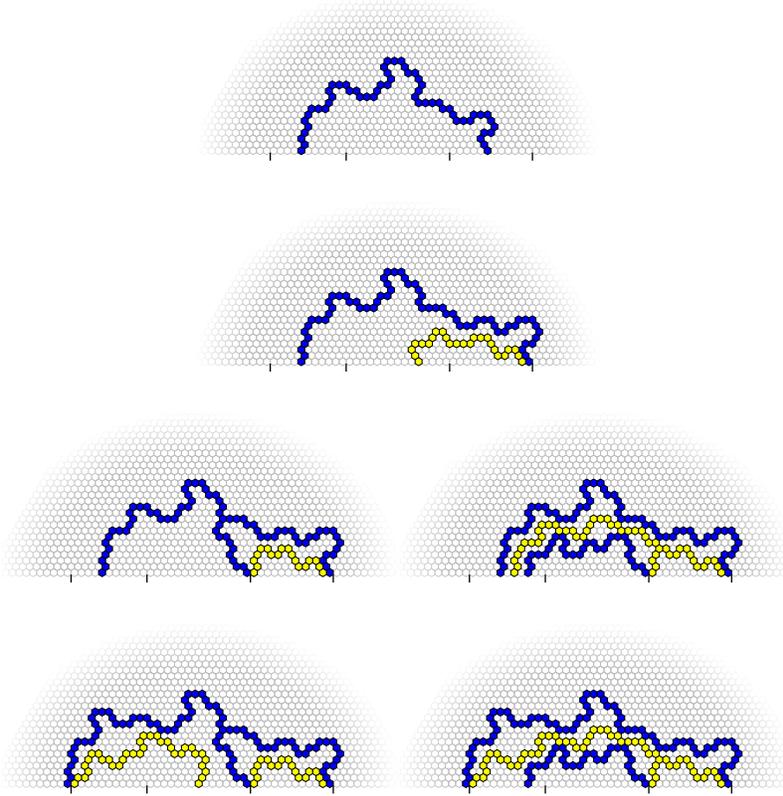}

\caption{The discrete triple partial derivative of the crossing
probability is sum of the probabilities of two multi-arm events.
The panels illustrate the events whose probability is
$P_C[x_1,x_2,x_3,x_4]$ (top), $\Delta_{x_4} P_C[x_1,x_2,x_3,x_4]$
(second), $\Delta_{x_3} \Delta_{x_4} P_C[x_1,x_2,x_3,x_4]$ (third row),
and $-\Delta_{x_1} \Delta_{x_3} \Delta_{x_4} P_C[x_1,x_2,x_3,x_4]$
(bottom row).}
\label{fig:dddcardy}
\end{figure}

Consider percolation on a half-plane triangular lattice, as in the
previous subsection, and let $P_C[x_1,x_2,x_3,x_4]$ be the probability
of at least one blue cluster spanning the intervals $(x_1,x_2)$ and
$(x_3,x_4)$; see the upper image in Figure \ref{fig:dddcardy}. Then
%
%
$\Delta_{x_4} P_C[x_1,x_2,x_3,x_4] = P_C[x_1,x_2,x_3,x_4] -
P_C[x_1,x_2,x_3,x_4-1]$ gives the probability that there is a crossing
for $(x_1,x_2,x_3,x_4)$ but \textit{not} $(x_1,x_2,x_3,\allowbreak x_4-1)$.
This event is represented by the second image in Figure \ref{fig:dddcardy}.
Similarly,
\[
-\Delta_{x_1} \Delta_{x_3} \Delta_{x_4}P_C[x_1,x_2,x_3,x_4]
\]
gives the probability that
\textit{one of the two} multi-arm events in the bottom image in
Figure \ref{fig:dddtripod} occurs. The event of a crossing cluster is
equivalent
to the event that one of these multi-arm events occurs for
\textit{some} (necessarily unique) set of three points $(a,b,c) \in
(x_1,x_2] \times(x_2,x_3] \times(x_3,x_4]$:
$a$ is (just right of)
the rightmost boundary point of the crossing cluster(s) in the
interval $(x_1,x_2)$; $b$ is (just right of) the rightmost point of
the crossing cluster in $(x_2,x_3)$ [if it exists; otherwise $b$ is
the rightmost boundary point in $(x_2,x_3)$ of a crossing yellow
cluster, as shown]; and $c$~is (just right of) the leftmost point of
the cluster(s) in $(x_3,x_4)$.

Thus $-\Delta_{x_1} \Delta_{x_3}
\Delta_{x_4}P_C[x_1,x_2,x_3,x_4]$ decomposes into the probabilities of two
multiarm events, the first of which is $-\Delta_{x_1} \Delta_{x_3}
\Delta_{x_4}P_T[x_1,x_2,x_3,x_4]$.

\subsection{Multi-arm events and the interface} \label{ss:interface}

\begin{figure}

\includegraphics{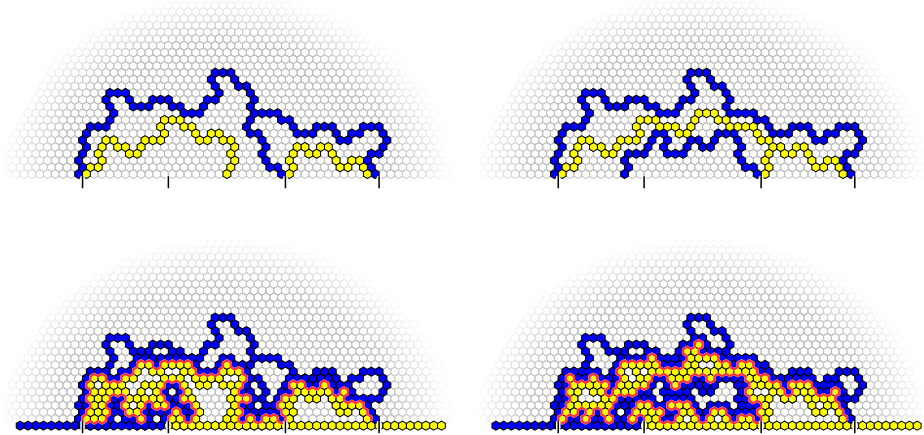}

\caption{The interface interpretation of the multi-arm events.}
\label{fig:interface}
\end{figure}

Consider the setting of Figures~\ref{fig:dddtripod}
and~\ref{fig:dddcardy}, and suppose we add an additional boundary
layer\vadjust{\goodbreak} of blue hexagons to the left of $x_2$ and yellow hexagons to the
right of $x_2$. Then let~$\gamma_{\mathrm{discrete}}$ be the discrete
interface starting at $x_2$. (See Figure \ref{fig:interface}.)

%
%

Then the union of the two multi-arm events at the bottom of
Figure \ref{fig:dddcardy}~de\-scribes the event that that $c$ is the first
boundary point that $\gamma_{\mathrm{discrete}}$ hits outside the
interval $(x_1,x_3)$, and that $a$ and $b$ are the leftmost and
rightmost boundary points hit by $\gamma_{\mathrm{discrete}}$
\textit{before} $c$. The left figure corresponds to~the~ca\-se that
$a$ is hit before $b$, and the right figure to the case that $b$ is
hit before~$a$.

\subsection{Continuum Watts' formula: A statement about SLE} \label
{ss.SLEstatement}

Like Cardy's formula, Watts' formula has a continuum analog, which is
a statement strictly about SLE$_6$. Fix real numbers $x_1 < x_2 < x_3
< x_4$, and let $s$ be their cross ratio. Consider the usual SLE$_6$
in the upper half-plane, where the starting point of the path is
$x_2$. Before Smirnov proved Cardy's formula for the scaling limit of
triangular lattice percolation, it was already known by Schramm that
$\cardy(s)$ represents the probability that $\gamma$ hits $(x_3,x_4)$
\textit{before} hitting $\R\setminus[x_1, x_4]$. [In the discrete
setting of Section \ref{ss:interface}, having $\gamma_{\mathrm{discrete}}$
hit $(x_3,x_4)$ before the complement of $(x_1,x_4)$ is equivalent to
the existence of a crossing.] In light of Section \ref{ss:interface}, the
following is the natural continuum analog of the tripod formula.
\begin{theorem} \label{t:continuumtripod}
Let $\SLEtripod(s)$ be the probability that both: 
%
\begin{longlist}[(2)]
\item[(1)] $\gamma$ first hits $(x_3,x_4)$ (at some time $t$) \textit
{before} it first hits $\R\setminus(x_1,x_4)$, and
\item[(2)] $\gamma$ hits the leftmost point of $\R\cap\gamma[0,t)$
before it hits the rightmost point.
\end{longlist}
Then $\SLEtripod(s) = \tripod(s)$.\vadjust{\eject}
\end{theorem}

Theorem \ref{t:continuumtripod} is the actual statement that was
proved by
Dub\'{e}dat, and the statement whose proof was sketched by Oded.
Dub\'{e}dat claimed further that Theorem~\ref{t:continuumtripod} would imply
the tripod formula (and hence Watts' formula) for the scaling limit of
critical triangular lattice percolation if one used the (at the time
\mbox{unpublished}) proof that SLE$_6$ is the scaling limit of the interface
[\citet{DubedatWatts}]. To be fully precise, one needs slightly more
than the fact that the interface scaling limit is SLE$_6$: it is
important to know that the discrete interface is unlikely to get close
to the boundary without hitting it. [Similar issues arise when using
Cardy's formula to prove SLE$_6$ convergence; see, e.g.,
\citet{CM}.]
Rather than address this (relatively minor technical) point here, we
will proceed to prove Theorem \ref{t:continuumtripod} in the manner outlined
by Oded and defer this issue until Section \ref{s.percolationstatement}.

It is convenient to have a name for the SLE versions of the multi-arm
events in Figure \ref{fig:dddcardy}. Say that a triple of distinct
real numbers $(a,b,c)$ with $a<0<b$
constitutes a \textit{tripod set} for $\gamma$ if for some $t > 0$ we
have:
\begin{longlist}[(2)]
\item[(1)] $\gamma(t) = c$;
\item[(2)] $\inf(\gamma[0,t)\cap\R) = a$;
\item[(3)] $\sup(\gamma[0,t)\cap\R) = b$.
\end{longlist}
There are a.s. a countably infinite number of tripod sets, but if
$x_1<0$ and $x_3 > 0$ is fixed, there is a.s. exactly one for which
$x_1<a<0<b<x_3$ and $c\notin(x_1,x_3)$. Let $\mathcal P_C\dvtx \R^3\to\R$
be the probability density function for this $(a,b,c)$.
(We see in Lemma \ref{density} that this density function exists.)
There are
also two types of tripod sets $(a,b,c)$: those for which $\gamma$ hits
$a$ first and those for which $\gamma$ hits $b$ first. Write
$\mathcal P_C = \mathcal P_A + \mathcal P_B$, where $\mathcal P_A$ and
$\mathcal P_B$ are the corresponding probability densities for
$a$-first and $b$-first tripod sets.

Now, we claim the following:
\begin{lemma}\label{density}
Using the notation above, the density functions $\mathcal P_C$ and
$\mathcal P_A$ exist, and
\[
\cardy(s) = \int_{x_1}^{0}\int_{0}^{x_3} \int_{x_3}^{x_4} \mathcal
P_C(a,b,c) \,dc \,db \,da
\]
and
\[
\SLEtripod(s) = \int_{x_1}^{0}\int_{0}^{x_3} \int_{x_3}^{x_4}
\mathcal P_A(a,b,c) \,dc \,db \,da.
\]
\end{lemma}
\begin{pf}
In the event that the density functions do not exist, we abuse
notation and let ``$\mathcal P_C(a,b,c) \,dc \,db \,da$'' and
``$\mathcal P_A(a,b,c) \,dc \,db \,da$'' denote the relevant measures,
which must exist. It is easy to see that the event that $(x_3,x_4)$
is hit before $\R\setminus[x_1,x_4]$ is equivalent to the event
that $(a,b,c) \in(x_1,0) \times(0,x_3) \times(x_3,x_4)$. By
definition, $\SLEtripod(s)$ is the probability of the same event
intersected with the event that $\gamma$ hits $a$ first. Finally,
observe that Cardy's formula of the cross-ratio of 4-points is
three-times differentiable, so the density function $\mathcal
P_C(a,b,c)$ exists and, consequently, the density function $\mathcal
P_A(a,b,c)$ also exists.
\end{pf}

Of course, from this, one has the immediate corollary:
\begin{corollary} \label{c.triplederivatives}
Using the above notation,
\[
\partial_{x_1} \partial_{x_3} \partial_{x_4} \cardy( \cratio
(x_1,0,x_3,x_4)) = \mathcal P_C(x_1,x_3,x_4)
\]
and
\[
\partial_{x_1} \partial_{x_3} \partial_{x_4} \SLEtripod(\cratio
(x_1,0,x_3,x_4)) = \mathcal P_A(x_1,x_3,x_4).
\]
\end{corollary}

If we could show further that
%
%
\begin{equation}\label{PA} \mathcal P_A(x_1,x_3,x_4) = \partial_{x_1}
\partial_{x_3} \partial_{x_4} \tripod(\cratio(x_1,0,x_3,x_4)),
\end{equation}
then this corollary and standard integration would imply
Theorem \ref{t:continuumtripod}, since we know that $\tripod(\cratio
(\cdot)) =
\SLEtripod(\cratio(\cdot))$ on the bounding planes \mbox{$x_1=0$}, \mbox{$x_2=0$}
and $x_3 = x_4$. Since we already have an explicit formula for
$\tripod$, the only remaining step is to explicitly compute $\mathcal P_A$.
Oded's approach is to compute the ratio $\mathcal P_A/\mathcal P_C$ as
the conditional probability [given that $(a,b,c)$ form a tripod set]
that $\gamma$ hits $a$ before $b$. Since $\mathcal P_C$ is known,
this determines~$\mathcal P_A$.

\section{Conditional probability that $a$ is hit first} \label{sec:cond-prob}

Schramm was very adept with using Mathematica to calculate all manner
of things. He probably would have considered this last step to be
routine, since it was for him straightforward to set up the right
equations and then let Mathematica solve them. At this point we refer
to his original Mathematica notebook from 2004, and explain the
various steps in the calculation. To be consistent with Oded's
notation, we now make the following substitutions:
\[
v_3 = a,  \qquad   W = x_2,  \qquad   v_1 = b,  \qquad    v_2 = c.
\]
(We assume $v_3 < W < v_1 < v_2$. Oded apparently chose this notation
because under cyclic reordering it was the same as $W=v_0,v_1,v_2,v_3$.)

First we formally define the function $\cardy(s)$ as in
(\ref{e.cardy}). The Mathematica function \Mcardy defined here
involves an additional parameter $\kappa$, but it specializes to
$\cardy(s)$ when $\kappa=6$. This more general formula is analogous
to Cardy's formula but gives the (conjectural) crossing probability
for the critical Fortuin--Kasteleyn random cluster model [with
$q=4\cos^2(4\pi/\kappa)$] with alternating wired-free-wired-free
boundary conditions. (See Rohde and Schramm [(\citeyear{RS}), conjecture 9.7],
for some background.)
This formula was known at Microsoft in 2003,
and most likely Oded copied it from another Mathematica notebook.
This formula was later independently discovered [\citet{BBK}]
(nonrigorously) and [\citet{DEuler}] (rigorously).\vspace*{6pt}

\begin{center}

\includegraphics{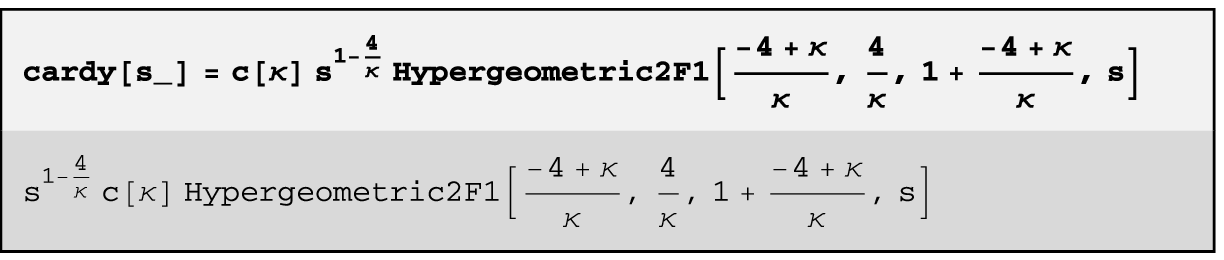}
\vspace*{2pt}
\end{center}

Consider the evolution of chordal SLE$_6$ started from $W$ and run to
$\infty$, when at time zero there are 3 marked points at positions
$v_1$, $v_2$ and $v_3$. We then let $W(t)$ represent the SLE$_6$
driving function [i.e., $W(t) = W(0) + \sqrt6 B_t$ where $B_t$ is a
standard Brownian motion] of the Loewner evolution
\[
\partial_t g_t(z) = \frac{2}{g_t(z)-W(t)},
\]
and interpret the $v_i$ as functions of $t$, evolving under the
Loewner flow, that is, $v_i(t):= g_t(v_i(0))$.

If $f$ is any function of $v_1, v_2, v_3, W$, we define
%
%
\begin{equation} \label{e.generator}
L(f):=  \frac{\partial}{\partial t} \mathbb E [ f(W(t),
v_1(t), v_2(t), v_3(t))] \bigg|_{t=0}.
\end{equation}
This is a new function of the same four variables which can be
calculated explicitly using It\^{o}'s formula as
\[
L(f) = \frac{\kappa}{2} \,\frac{\partial^2}{\partial W^2}f + \sum
_{i=1}^3 \frac{2 ({\partial}/{\partial v_i})f}{v_i - W}.
\]
This operator is defined as \ML in the Mathematica code below.\vspace*{6pt}
\begin{center}

\includegraphics{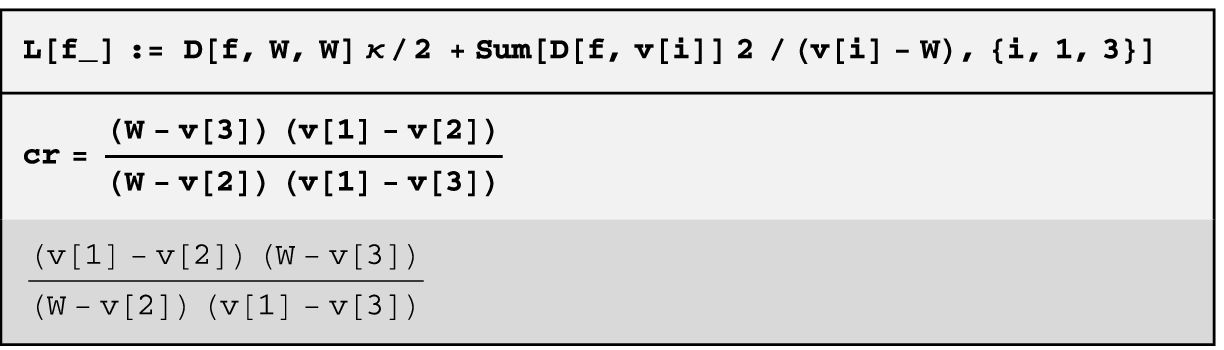}
\vspace*{2pt}
\end{center}
Similarly in the Mathematica code, \Mcr is the cross-ratio [defined in
(\ref{e.cratio})], that is,
%
%
\begin{equation}\label{e.cratio} \cratio:=\cratio(v_3,W, v_1, v_2)=
\frac{(W-v_3)(v_1-v_2)}{(W - v_2)(v_1-v_3)}.
\end{equation}

In the next line, Oded performed a consistency check.
Cardy's formula should be a martingale for the SLE$_\kappa$ diffusion,
hence $L(\cardy(\cratio(\cdot))) = 0$.\vspace*{3pt}
\begin{center}

\includegraphics{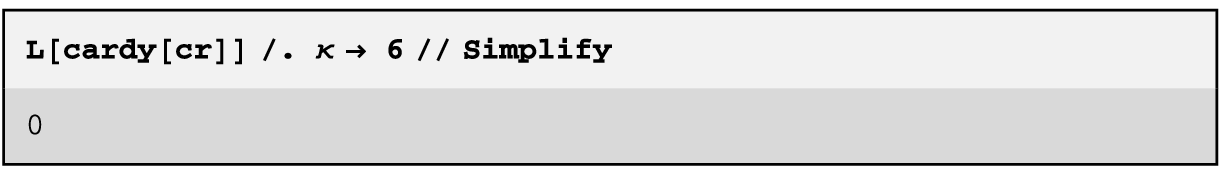}
\vspace*{3pt}
\end{center}
Next, Oded computes the triple derivative \Mh of Cardy's formula
(which is the same as the $\mathcal P_C$ defined in
Section \ref{ss.SLEstatement}).\vspace*{3pt}

\begin{center}

\includegraphics{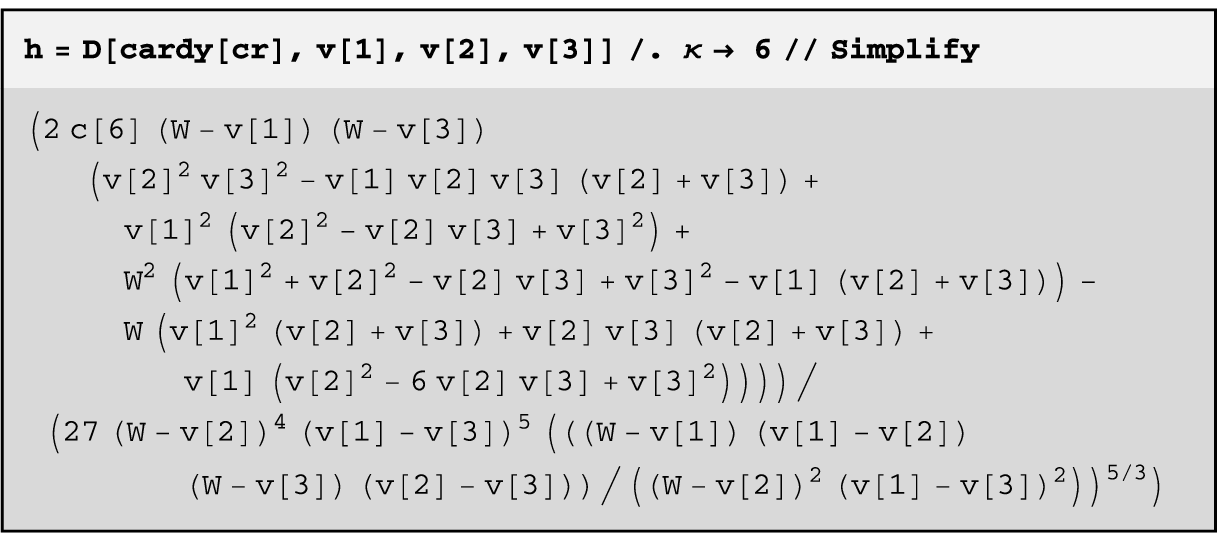}
\vspace*{3pt}
\end{center}

That is, he computes
\[
h(v_3,W,v_1,v_2) := \partial_{v_1} \partial_{v_2} \partial_{v_3}
\cardy\biggl(\frac{(W-v_3)(v_1-v_2)}{(W-v_2)(v_1-v_3)}\biggr).
\]
The result is somewhat complicated, but we may ignore it, since it is
just an
intermediate result.

The next step involves conditioning on an event of zero probability,
the event that $(v_3,v_1,v_2)$ is a tripod set. We can make sense of
this by introducing a triple difference of Cardy's formula and
recalling the results of Section \ref{ss.SLEstatement}. First we introduce
notation to describe some small evolving intervals. For given values
$v_1(0), v_2(0), v_3(0), W(0)$, pick $\eps$ small enough so that the
intervals $(v_i(0), v_i(0)+\eps)$ are disjoint and do not contain
$W(0)$. Write $\tilde v_i(0) = v_i(0) + \eps$. Define $\tilde
v_i(t)$ using the Loewner evolution, and write $\eps_i(t) := \tilde
v_i(t) - v_i(t)$. Let us write
%
%
\begin{eqnarray}\label{heps}
&&h_{\eps_1, \eps_2, \eps_3}(v_3,W,v_1,v_2)\nonumber\\[-8pt]\\[-8pt]
&&\qquad := \Delta_{v_1}^{(\eps
_1)} \Delta_{v_2}^{(\eps_2)} \Delta_{v_3}^{(\eps_3)} \cardy
(\cratio(v_3,W,v_1,v_2)),\nonumber
\end{eqnarray}
where $\Delta_v^{(\eps)}$ is the difference operator defined by
\[
\Delta_v^{(\eps)} f(v) = f(v+\eps)-f(v).
\]
Note that the $\Delta^{(\eps_i)}_{v_i}$ depend on $t$. By Corollary
\ref{c.triplederivatives}, equation (\ref{heps}) at time~$t$
represents the conditional probability (given the Loewner evolution up
to time $t$) that there is a tripod set in $[v_3(0), \tilde v_3(0)]
\times[v_1(0), \tilde v_1(0)] \times[v_2(0), \tilde v_2(0)]$.
By Girsanov's theorem,
conditioning on this event induces a~drift on the Brownian motion $W_t$
driving the SLE, where the drift is
\[
\kappa  \,\partial_W \log h_{\eps_1, \eps_2, \eps_3}(v_3,W,v_1,v_2).
\]
Observe that
\begin{eqnarray*}
\frac{\partial}{\partial W} \log h_{\eps_1, \eps_2, \eps_3} &=&
\frac{({\partial}/{\partial W}) h_{\eps_1,\eps_2,\eps
_3}}{h_{\eps_1,\eps_2,\eps_3}} =
\frac{\iiint({\partial}/{\partial W}) h}{\iiint h} \\
&\approx&\frac
{\eps_1 \eps_2 \eps_3 ({\partial}/{\partial W}) h}{\eps_1\eps
_2\eps_3 h} = \partial_W \log h,
\end{eqnarray*}
where there triple integral is over $\prod[v_i, \tilde v_i]$.
Thus upon taking the limit $\eps\to0$, the drift becomes
%
%
\begin{equation} \label{e.drift}\drift(t):= \kappa  \partial_W
\log h(v_3,W,v_1,v_2).
\end{equation}
The next Mathematica code explicitly computes (\ref{e.drift}).\vspace*{3pt}

\begin{center}

\includegraphics{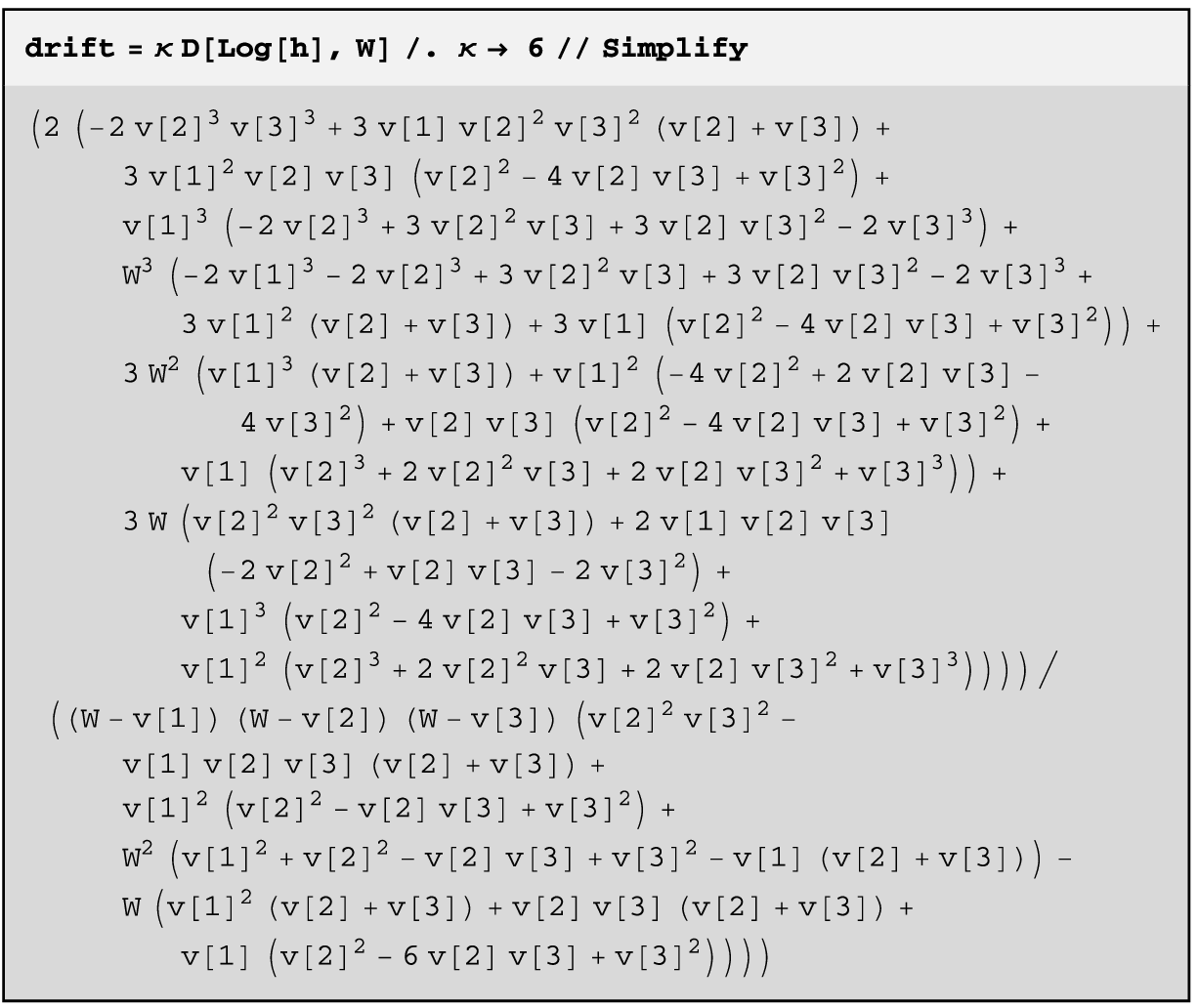}

\end{center}

The expression above is complicated, but again it is an intermediate
result that we do not need to calculate or read ourselves. The first
line of the Mathematica code below defines the generator \MLd(which
we will write as~$L_1$) for the conditioned SLE$_6$, where the driving
function $W_t$ has the drift given above. Here $L_1$ is defined as in
(\ref{e.generator}) except that the expectation is with respect to the
law of $W_t$ with the drift term (\ref{e.drift}). Thus
\[
L_1(f) :=
L(f) + \drift(t)  \,\partial_W f.
\]
As before, if $f$ is a real function
of $W,v_1,v_2,v_3$, then $L_1(f)$ will be a function of the same four
variables.

We now compute the probability in the modified diffusion that
$v_3$ is absorbed before $v_1$, that is, that $W(t)$ collides with $v_3(t)$
before colliding with~$v_1(t)$. This probability will be a martingale
that only depends upon the cross-ratio $s$. Thus, in the next
paragraph, we specialize and consider functions of~$W,v_1,v_2,v_3$
that have the form $f(\cratio(v_3,W,v_1,v_2))$ where $f\dvtx\R\to\R$ is a
function of \textit{one} variable. We would like to find a one-parameter
function~$f$ for which $f(\cratio(v_3,W,v_1,v_2))$ is a
martingale with respect to this modified diffusion, so we will require
that $L_1(f(\cratio(v_3,W,v_1,v_2))) = 0$. What one-parameter
functions $f$
have this property?

Oded answers this question with some clever Mathematica work. First,
he re-expresses the differential equation $L_1(f(\cratio
(v_3,W,v_1,v_2))) = 0$---which involves the four parameters
$W,v_1,v_2,v_3$---in terms of the
parameters $s,v_1,v_2,\allowbreak v_3$. He does this by setting $s$ equal to the
expression for $\cratio$ given in (\ref{e.cratio}), solving to get $W$
in terms of the other variables and plugging this new expression for
$W$ into the expression $L_1(f(\cratio(v_3,W,v_1,v_2)))$.\vspace*{6pt}

\begin{center}

\includegraphics{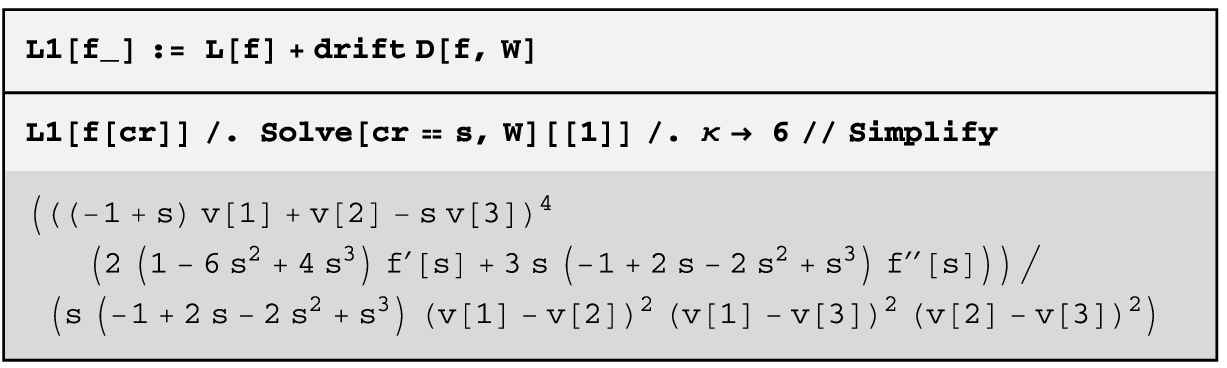}
\vspace*{6pt}
\end{center}

This expression for $L_1(f(\cratio(v_3,W,v_1,v_2)))$ depends on
$f'$ and $f''$, and equating it to zero yields a differential
equation for $f$, the unknown one-parameter function of the
cross-ratio that we seek,
\[
2(1-6s^2+4s^3)f'(s) + 3 s (-1+2s-2s^2+s^3) f''(s) = 0.
\]
Oded solves this differential equation, which yields the function
$f$ up to two free parameters $C_1$ and $C_2$.

\begin{center}

\includegraphics{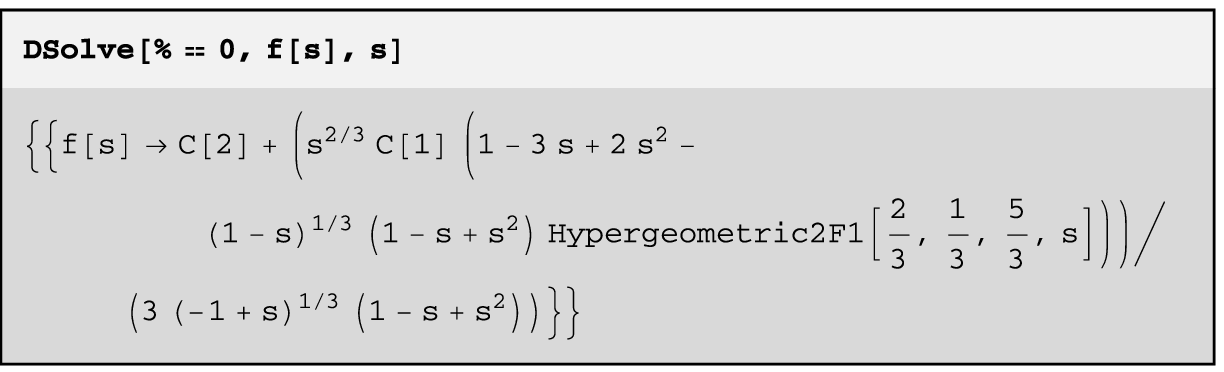}
\vspace*{3pt}
\end{center}

Here Mathematica gives
\[
f(s) = C_2 + C_1 s^{2/3} \frac{1-3 s+2 s^2-(1-s)^{1/3} (1-s+s^2)
\,\hspace*{0.4pt}{}_2\hspace*{-0.4pt}F_1(2/3,1/3;5/3;s)}{3 (-1+s)^{1/3} (1-s+s^2)}.
\]

The conditional probability that we seek tends to $1$ when $s\to0$ and
tends to $0$
when $s\to1$, and this determines $C_1$ and $C_2$: $C_2$ must be $1$,
and~$C_1$ follows from Gauss's hypergeometric formula,
\[
{}_2\hspace*{-0.4pt}F_1(a,b;c;1) = \frac{\Gamma(c)\Gamma(c-a-b)}{\Gamma(c-a)\Gamma
(c-b)}.
\]
Solving for $C_1$ and substituting, we find that the conditional
probability that $\gamma$ hits $v_3$ before $v_1$, given that
$(v_3,v_1,v_2)$ is
a tripod set, is given by
\[
f(s) = 1 - \frac{\Gamma(4/3)}{\Gamma(2/3)\Gamma(
5/3)}   s^{2/3}  \biggl[\frac{-1+3 s-2 s^2}{(1-s)^{1/3} (1-s+s^2)} +
{}_2\hspace*{-0.4pt}F_1\biggl(\frac23,\frac13;\frac53;s\biggr)\biggr],
\]
where $s$ is the cross-ratio of $v_3,0,v_1,v_2$.

\section{Comparison of triple derivatives}

Taking $\mathcal P_A$ and $\mathcal P_C$ as defined in
Section \ref{ss.SLEstatement}, and $f$ and $h$ as defined in the previous
section, we now have
\[
h(v_3,0,v_1,v_2) = \mathcal P_C(v_3,v_1,v_2)
\]
and
\[
f(\cratio(v_3,0,v_1,v_2))=\mathcal P_A(v_3,v_1,v_2)/\mathcal P_C(v_3,v_1,v_2).
\]
In principle the next step toward proving (\ref{PA}) (and hence
Theorem \ref{t:continuumtripod}) would be to integrate $\mathcal P_A =
f(\cratio(\cdot))h$ over the three variables $v_1,v_2,v_3$ and show
that one obtains $\tripod(\cratio(\cdot))$. In Oded's original notes,
he stated that this could be done using integration by parts.
Fortunately (for those who lack Oded's skill at integrating) we
already know (thanks to Watts) what we expect $\tripod$ to be, so we
can instead differentiate $\tripod(\cratio(\cdot))$ three times
(w.r.t. $v_1,v_2,v_3$), and check that it equals
$f(\cratio(\cdot))h$. The Mathematica code in this final section was
generated by the authors of this paper, not by Schramm.

First we redefine \Mcardy to have an explicit constant and define the
purported tripod probability.

\begin{center}

\includegraphics{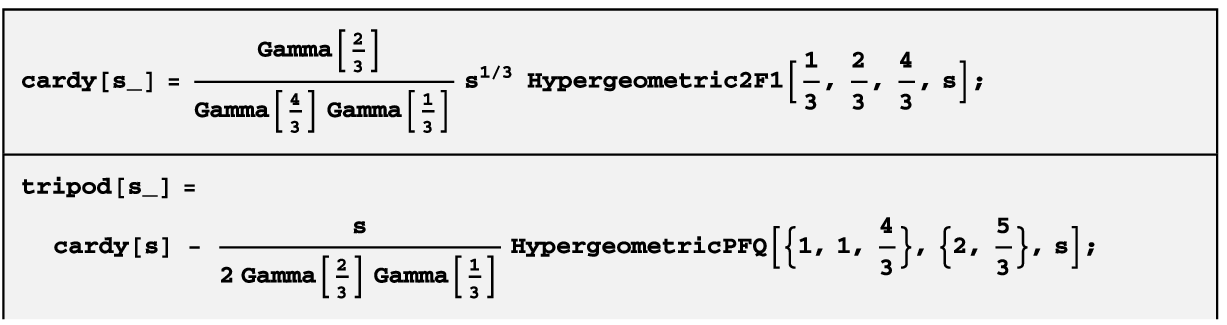}
\vspace*{3pt}
\end{center}

Next we differentiate Cardy's formula three times.\vspace*{6pt}

\begin{center}

\includegraphics{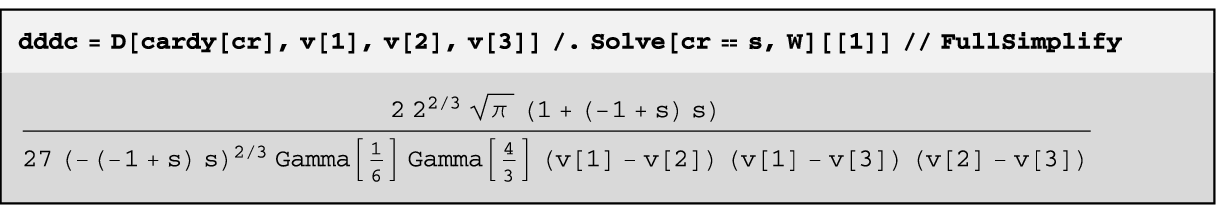}
\vspace*{3pt}

\end{center}

This is the same as \Mh defined earlier, but with the trick of
eliminating the variable $W$ and expressing the formula in terms of $s$.
Next we triply differentiate the purported tripod probability.\vspace*{6pt}

\begin{center}

\includegraphics{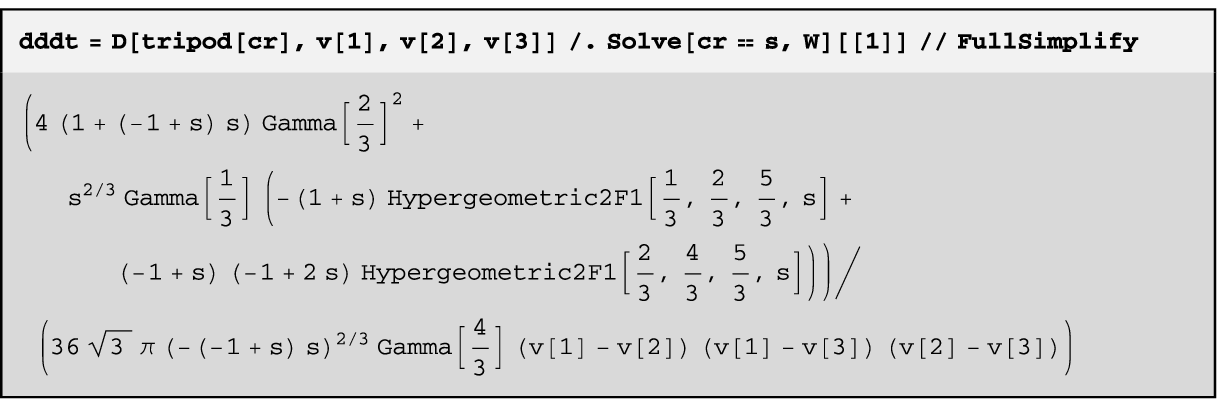}
\vspace*{3pt}

\end{center}

Notice that the triple derivative of the tripod probability is
expressed in terms of two different hypergeometric functions. In
order to compare this expression with the conditional probability
computed in Section \ref{sec:cond-prob}, we need to use some
hypergeometric identities.
We use one of Gauss's relations between ``contiguous''
hypergeometric functions [Erd{\'e}lyi et al.
(\citeyear{EMOT}), Section~2.8, equation 33], to write
\[
-\tfrac13 \,\hspace*{0.4pt}{}_2\hspace*{-0.4pt}F_1\bigl(\tfrac23,\tfrac43;\tfrac53;s\bigr) + \tfrac23 (1-s)
\,\hspace*{0.4pt}{}_2\hspace*{-0.4pt}F_1\bigl(\tfrac53,\tfrac43;\tfrac53;s\bigr) - \tfrac13
\,\hspace*{0.4pt}{}_2\hspace*{-0.4pt}F_1\bigl(\tfrac23,\tfrac
13;\tfrac53;s\bigr)=0.
\]
But ${}_2\hspace*{-0.4pt}F_1(c,b;c;s) = (1-s)^{-b}$ \citet{EMOT}, Section 2.8,
equation 4, so
%
%
\begin{equation}
{}_2\hspace*{-0.4pt}F_1\bigl(\tfrac23,\tfrac43;\tfrac53;s\bigr) = 2 (1-s)^{-1/3}
- {}_2\hspace*{-0.4pt}F_1\bigl(\tfrac
23,\tfrac13;\tfrac53;s\bigr).
\end{equation}

\begin{center}

\includegraphics{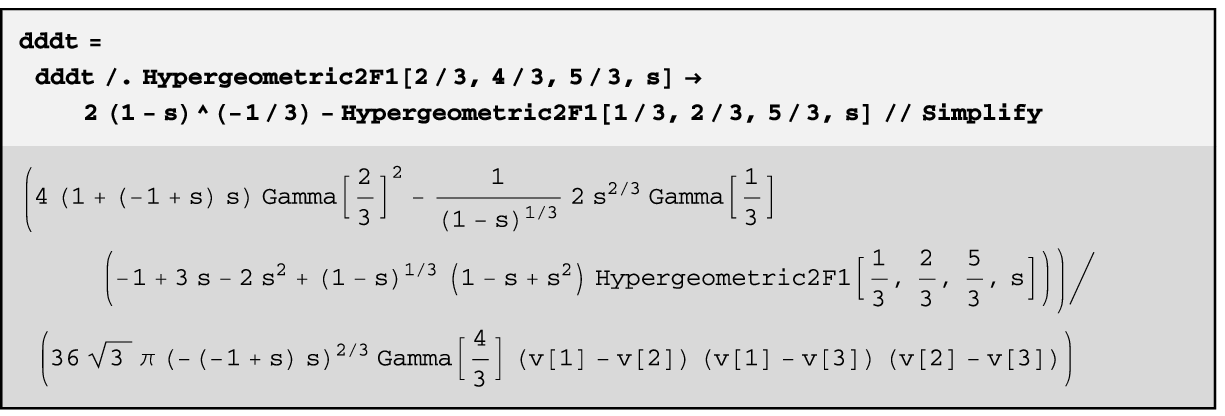}
\vspace*{3pt}

\end{center}

Next we compare the two expressions for the conditional probability and
verify that they are the same.\vspace*{6pt}

\begin{center}

\includegraphics{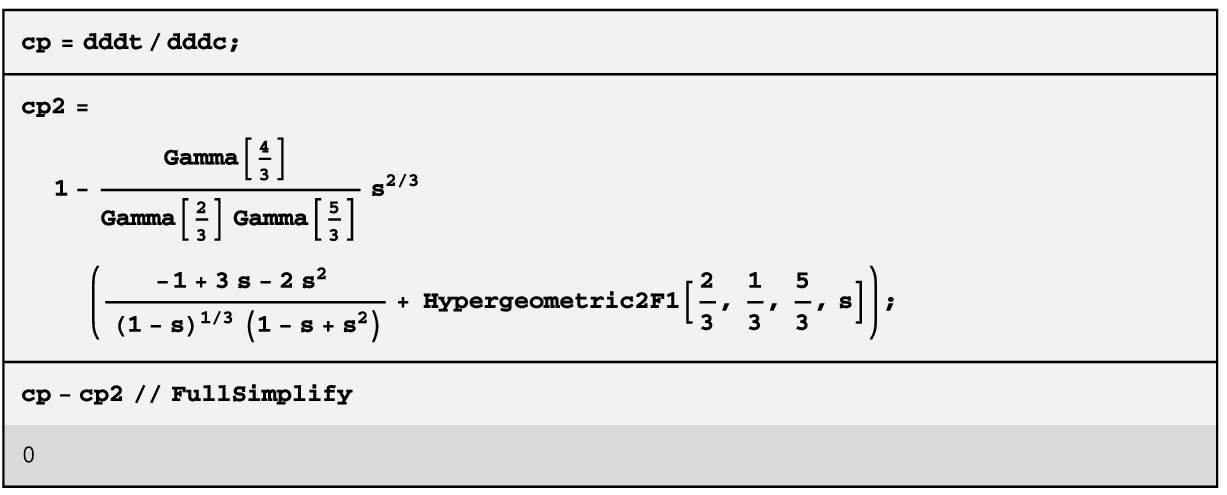}
\vspace*{3pt}

\end{center}

Therefore the triple derivatives agree, and we have established (\ref{PA}).

\section{Percolation statement} \label{s.percolationstatement}

\mbox{\hspace*{8pt}}We have the established equivalence of\break $\SLEtripod(s)$ and $\tripod
(s)$, but we still need
to make the connection to percolation.
\begin{theorem} \label{t.watts}
Let $D \subset\C$ be a fixed bounded
Jordan domain with marked points $x_1,x_2,x_3,x_4$ on its boundary.
For any $\eps$, we may consider the hexagonal lattice rescaled to
have side length $\eps$ and color the faces blue and yellow
according to site percolation. Let $B$ be the closure of the set of
blue faces, and let $P^\eps$ be the probability that $B \cap\overline D$
contains a connected component that intersects all four boundary
segments $(x_1,x_2)$, $(x_2,x_3)$, $(x_3,x_4)$ and $(x_4,x_1)$.
Then
\[
\lim_{\eps\to0} P^\eps= \watts(s).
\]
\end{theorem}

Proving Theorem \ref{t.watts} solves the problem addressed by Watts.
However, we remark that more general statements are probably possible.
Any domain~$D$ with four marked boundary points has a ``center''
$c(D)$ with the property that a conformal map taking the domain to a
rectangle (and the points to the corners) sends $c$ to the center of
the rectangle. Oded would probably have preferred to show that for any
sequence $D_n$ of simply connected \textit{marked hexagonal domains}
(domains comprised of unions of hexagons within a~fixed hexagonal lattice $H$ with four
marked boundary points of cross ratios $s_n$ converging to $s$), the
probability of the Watts event tends to $\watts(s)$ provided that the
distance from $c(D_n)$ to $\partial D_n$ tends to~$\infty$. (Oded's
SLE convergence results are similarly general
[\citet{LSW}, Schramm and Sheffield
(\citeyear{SchrammSheffieldHE,SchrammSheffieldGFF})].) However, Oded's
derivation of Watts' formula (like Dub\'{e}dat's derivation) depends on
SLE$_6$ convergence, and existing SLE$_6$ convergence statements [e.g.,
\citet{CM}] are not quite general enough to imply this.
\begin{pf*}{Proof of Theorem \ref{t.watts}}
As shown in Section \ref{s.tripodreduction}, it suffices to prove the
analogous statement about tripod events $(x_1,x_2)$, $(x_2,x_3)$,
and $(x_3,x_4)$ and the function $\tripod(s)$.

Let $\eps_n$ be a sequence of positive reals tending to zero, and
define $\gamma_n$ to be the random interface in $D$ obtained from
percolation on $\eps_n$ times the hexagonal lattice, between the
lattice points closest to $x_2$ and $x_4$. Let~$a_n,\allowbreak b_n,c_n$ be the
tripod set for this interface and the points $x_1$ and $x_3$, that
is,~$c_n$ is the first point on $\gamma_n \cap\partial D_n$
outside the boundary segment $(x_1,x_3)$, and the interface
$\gamma_n$ up to point $c_n$ last hits the boundary intervals
$(x_1,x_2)$ and $(x_2,x_3)$ at $a_n$ and $b_n$, respectively. From the work
of \citet{CM}, we can couple the $\gamma_n$ and
$\gamma$ in such a way that $\gamma_n \to\gamma$ almost surely in
the uniform topology (in which two curves are close if they can be
parameterized in such a way that they are close at all times). By
the compactness of $\partial D$ (and the corresponding compactness---in
the topology of convergence in law---of the space of
measures on $\partial D$) it is not hard to see that there must be a
subsequence of the $n$ values and a coupling of the $\gamma_n$
with~$\gamma$ in which the entire quadruple $(\gamma_n, a_n, b_n, c_n)$
converges almost surely to some limit. If we could show further
that this limit must be $(\gamma, a,b,c)$ almost surely, this would
imply the theorem, since uniform topology convergence would imply
that if $\gamma$ hits $a$ before $b$ then $\gamma_n$ hits $a_n$
before $b_n$ for large enough $n$ almost surely. However, it is not
clear a priori that this limit is $(\gamma, a,b,c)$
almost surely (even though the $\gamma_n$ converge to $\gamma$),
since while $\gamma$ touches the boundary at $a$, $b$ and $c$, it
could be that $\gamma_n$ comes close to the boundary at these points
without touching it.

To obtain a contradiction, let us suppose that there is a uniformly positive
probability (i.e., bounded away from $0$ as $n\to\infty$) that,
say, the limit of the~$a_n$ is not $a$. (The
argument for the $b_n$ and the $c_n$ is essentially the same.) Then
there must be an open interval $(\alpha_1, \alpha_2)$ of the
boundary and an open subinterval
$(\beta_1,\beta_2)\subset(\alpha_1,\alpha_2)$ of the boundary (with
$\beta_1\neq\alpha_1$ and $\beta_2\neq\alpha_2$) such that there is
a uniformly positive probability that $a$ lies $(\beta_1,\beta_2)$
but the limit of the $a_n$ does not lie in that\vadjust{\eject}
$(\alpha_1,\alpha_2)$. Now we can expand the Jordan domain $D$ to a
larger Jordan domain $\tilde D$ that includes a neighborhood of
$(\beta_1,\beta_2)$, but where the boundary of $\tilde D$ agrees
with boundary of $D$ outside of $(\alpha_1,\alpha_2)$. Let $\tilde
\gamma_n$ denote the discrete interfaces in this expanded domain.
We can couple the $\tilde\gamma_n$ with the $\gamma_n$ in such a
way that the two agree whenever~$\tilde\gamma_n$ does not leave $D$ (by
using the same percolation to define both). But now we have a
coupling of the $\tilde\gamma_n$ sequence with the property that
there is a~positive probability that the limit of the $\tilde
\gamma_n$ is a path that\vspace*{2pt} hits the boundary of $\tilde D \setminus D$
without entering $\tilde D \setminus D$. This implies that if the
$\tilde\gamma_n$ converge in law, they must converge to a random
path that\vspace*{2pt} with positive probability hits the boundary of $\tilde D
\setminus D$ without entering $\tilde D \setminus D$. By the
Camia--Newman theorem, applied to the domain $\tilde D$, the
$\tilde\gamma_n$ converge in law to chordal SLE$_6$ in $\tilde D$,
and on the event that SLE$_6$ hits $\partial(\tilde D \setminus D)$,
it will a.s. enter $\tilde D \setminus D$, a~contradiction.
\end{pf*}

\section*{Acknowledgments}
We thank Nike Sun, Kalle Kyt\"{o}l\"{a} and Stas Smirnov for their
comments on an earlier version of this article.


%

%
\printaddresses


\begin{thebibliography}{18}

\bibitem[\protect\citeauthoryear{Bauer, Bernard and Kyt{\"o}l{\"a}}{2005}]{BBK}
%
\begin{barticle}[mr]
\bauthor{\bsnm{Bauer},~\bfnm{Michel}\binits{M.}},
\bauthor{\bsnm{Bernard},~\bfnm{Denis}\binits{D.}} \AND
\bauthor{\bsnm{Kyt{\"o}l{\"a}},~\bfnm{Kalle}\binits{K.}}
(\byear{2005}).
\btitle{Multiple {S}chramm--{L}oewner \mbox{evolutions} and statistical mechanics
martingales}.
\bjournal{J. Stat. Phys.}
\bvolume{120}
\bpages{1125--1163}.
\bid{doi={10.1007/s10955-005-7002-5}, issn={0022-4715}, mr={2187598}}
\end{barticle}
%
\endbibitem

%
\bibitem[\protect\citeauthoryear{Camia and Newman}{2007}]{CM}
%
\begin{barticle}[mr]
\bauthor{\bsnm{Camia},~\bfnm{Federico}\binits{F.}} \AND
\bauthor{\bsnm{Newman},~\bfnm{Charles~M.}\binits{C.~M.}}
(\byear{2007}).
\btitle{Critical percolation exploration path and {${\rm SLE}\sb6$}: A~proof
of convergence}.
\bjournal{Probab. Theory Related Fields}
\bvolume{139}
\bpages{473--519}.
\bid{doi={10.1007/s00440-006-0049-7}, issn={0178-8051}, mr={2322705}}
\end{barticle}
%
\endbibitem

\bibitem[\protect\citeauthoryear{Cardy}{1992}]{Cardy}
%
\begin{barticle}[mr]
\bauthor{\bsnm{Cardy},~\bfnm{John~L.}\binits{J.~L.}}
(\byear{1992}).
\btitle{Critical percolation in finite geometries}.
\bjournal{J. Phys. A}
\bvolume{25}
\bpages{L201--L206}.
\bid{issn={0305-4470}, mr={1151081}}
\end{barticle}
%
\endbibitem

\bibitem[\protect\citeauthoryear{Dub{\'e}dat}{2004}]{dubedatBM}
%
\begin{barticle}[mr]
\bauthor{\bsnm{Dub{\'e}dat},~\bfnm{Julien}\binits{J.}}
(\byear{2004}).
\btitle{Reflected planar {B}rownian motions, intertwining relations and
crossing probabilities}.
\bjournal{Ann. Inst. Henri Poincar\'e Probab. Stat.}
\bvolume{40}
\bpages{539--552}.
\bid{doi={10.1016/j.anihpb.2003.11.005}, issn={0246-0203}, mr={2086013}}
\end{barticle}
%
\endbibitem

\bibitem[\protect\citeauthoryear{Dub{\'e}dat}{2006a}]{DubedatWatts}
%
\begin{barticle}[mr]
\bauthor{\bsnm{Dub{\'e}dat},~\bfnm{Julien}\binits{J.}}
(\byear{2006a}).
\btitle{Excursion decompositions for {SLE} and {W}atts' crossing formula}.
\bjournal{Probab. Theory Related Fields}
\bvolume{134}
\bpages{453--488}.
\bid{doi={10.1007/s00440-005-0446-3}, issn={0178-8051}, mr={2226888}}
\end{barticle}
%
\endbibitem

\bibitem[\protect\citeauthoryear{Dub{\'e}dat}{2006b}]{DEuler}
%
\begin{barticle}[mr]
\bauthor{\bsnm{Dub{\'e}dat},~\bfnm{Julien}\binits{J.}}
(\byear{2006b}).
\btitle{Euler integrals for commuting {SLE}s}.
\bjournal{J. Stat. Phys.}
\bvolume{123}
\bpages{1183--1218}.
\bid{doi={10.1007/s10955-006-9132-9}, issn={0022-4715}, mr={2253875}}
\end{barticle}
%
\endbibitem

\bibitem[\protect\citeauthoryear{Erd{\'e}lyi et al.}{1953}]{EMOT}
%
\begin{bbook}[auto:STB|2011-03-03|12:04:44]
\bauthor{\bsnm{Erd{\'e}lyi},~\bfnm{Arthur}\binits{A.}},
\bauthor{\bsnm{Magnus},~\bfnm{Wilhelm}\binits{W.}},
\bauthor{\bsnm{Oberhettinger},~\bfnm{Fritz}\binits{F.}} \AND
\bauthor{\bsnm{Tricomi},~\bfnm{Francesco~G.}\binits{F.~G.}}
(\byear{1953}).
\btitle{Higher Transcendental Functions}
\bvolume{I}.
\bpublisher{McGraw-Hill}, \baddress{New York}.
\bnote{Based, in part, on notes left by Harry Bateman.}
\end{bbook}
%
\endbibitem

\bibitem[\protect\citeauthoryear{Langlands et al.}{1992}]{LPPS}
%
\begin{barticle}[mr]
\bauthor{\bsnm{Langlands},~\bfnm{R.~P.}\binits{R.~P.}},
\bauthor{\bsnm{Pichet},~\bfnm{C.}\binits{C.}},
\bauthor{\bsnm{Pouliot},~\bfnm{Ph.}\binits{P.}} \AND
\bauthor{\bsnm{Saint-Aubin},~\bfnm{Y.}\binits{Y.}}
(\byear{1992}).
\btitle{On the universality of crossing probabilities in two-dimensional
percolation}.
\bjournal{J. Stat. Phys.}
\bvolume{67}
\bpages{553--574}.
\bid{issn={0022-4715}, mr={1171144}}
\end{barticle}
%
\endbibitem

\bibitem[\protect\citeauthoryear{Lawler, Schramm and Werner}{2004}]{LSW}
%
\begin{barticle}[mr]
\bauthor{\bsnm{Lawler},~\bfnm{Gregory~F.}\binits{G.~F.}},
\bauthor{\bsnm{Schramm},~\bfnm{Oded}\binits{O.}} \AND
\bauthor{\bsnm{Werner},~\bfnm{Wendelin}\binits{W.}}
(\byear{2004}).
\btitle{Conformal invariance of planar loop-erased random walks and uniform
spanning trees}.
\bjournal{Ann. Probab.}
\bvolume{32}
\bpages{939--995}.
\bid{doi={10.1214/aop/1079021469}, issn={0091-1798}, mr={2044671}}
\end{barticle}
%
\endbibitem

\bibitem[\protect\citeauthoryear{Maier}{2003}]{Maier}
%
\begin{barticle}[mr]
\bauthor{\bsnm{Maier},~\bfnm{Robert~S.}\binits{R.~S.}}
(\byear{2003}).
\btitle{On crossing event formulas in critical two-dimensional percolation}.
\bjournal{J.~Stat. Phys.}
\bvolume{111}
\bpages{1027--1048}.
\bid{doi={10.1023/A:1023006413433}, issn={0022-4715}, mr={1975920}}
\end{barticle}
%
\endbibitem

\bibitem[\protect\citeauthoryear{Rohde and Schramm}{2005}]{RS}
%
\begin{barticle}[mr]
\bauthor{\bsnm{Rohde},~\bfnm{Steffen}\binits{S.}} \AND
\bauthor{\bsnm{Schramm},~\bfnm{Oded}\binits{O.}}
(\byear{2005}).
\btitle{Basic properties of {SLE}}.
\bjournal{Ann. of Math. (2)}
\bvolume{161}
\bpages{883--924}.
\bid{doi={10.4007/annals.2005.161.883}, issn={0003-486X}, mr={2153402}}
\end{barticle}
%
\endbibitem

\bibitem[\protect\citeauthoryear{Schramm and Sheffield}{2005}]{SchrammSheffieldHE}
%
\begin{barticle}[mr]
\bauthor{\bsnm{Schramm},~\bfnm{Oded}\binits{O.}} \AND
\bauthor{\bsnm{Sheffield},~\bfnm{Scott}\binits{S.}}
(\byear{2005}).
\btitle{Harmonic explorer and its convergence to {${\rm SLE}\sb4$}}.
\bjournal{Ann. Probab.}
\bvolume{33}
\bpages{2127--2148}.
\bid{doi={10.1214/009117905000000477}, issn={0091-1798}, mr={2184093}}
\end{barticle}
%
\endbibitem

\bibitem[\protect\citeauthoryear{Schramm and Sheffield}{2009}]{SchrammSheffieldGFF}
%
\begin{barticle}[mr]
\bauthor{\bsnm{Schramm},~\bfnm{Oded}\binits{O.}} \AND
\bauthor{\bsnm{Sheffield},~\bfnm{Scott}\binits{S.}}
(\byear{2009}).
\btitle{Contour lines of the two-dimensional discrete {G}aussian free field}.
\bjournal{Acta Math.}
\bvolume{202}
\bpages{21--137}.
\bid{doi={10.1007/s11511-009-0034-y}, issn={0001-5962}, mr={2486487}}
\end{barticle}
%
\endbibitem

\bibitem[\protect\citeauthoryear{Simmons, Kleban and Ziff}{2007}]{SKZ}
%
\begin{barticle}[mr]
\bauthor{\bsnm{Simmons},~\bfnm{Jacob J.~H.}\binits{J.~J.~H.}},
\bauthor{\bsnm{Kleban},~\bfnm{Peter}\binits{P.}} \AND
\bauthor{\bsnm{Ziff},~\bfnm{Robert~M.}\binits{R.~M.}}
(\byear{2007}).
\btitle{Percolation crossing formulae and conformal field theory}.
\bjournal{J. Phys. A}
\bvolume{40}
\bpages{F771--F784}.
\bid{doi={10.1088/1751-8113/40/31/F03}, issn={1751-8113}, mr={2345284}}
\end{barticle}
%
\endbibitem

\bibitem[\protect\citeauthoryear{Smirnov}{2001}]{Smirnov}
%
\begin{barticle}[mr]
\bauthor{\bsnm{Smirnov},~\bfnm{Stanislav}\binits{S.}}
(\byear{2001}).
\btitle{Critical percolation in the plane: Conformal invariance, {C}ardy's
formula, scaling limits}.
\bjournal{C. R. Acad. Sci. Paris S\'er. I Math.}
\bvolume{333}
\bpages{239--244}.
\bid{doi={10.1016/S0764-4442(01)01991-7}, issn={0764-4442}, mr={1851632}}
\end{barticle}
%
\endbibitem

\bibitem[\protect\citeauthoryear{Watts}{1996}]{Watts}
%
\begin{barticle}[mr]
\bauthor{\bsnm{Watts},~\bfnm{G.~M.~T.}\binits{G.~M.~T.}}
(\byear{1996}).
\btitle{A crossing probability for critical percolation in two dimensions}.
\bjournal{J.~Phys. A}
\bvolume{29}
\bpages{L363--L368}.
\bid{doi={10.1088/0305-4470/29/14/002}, issn={0305-4470}, mr={1406907}}
\end{barticle}
%
\endbibitem

\bibitem[\protect\citeauthoryear{Zhan}{2008}]{Zhan}
%
\begin{barticle}[mr]
\bauthor{\bsnm{Zhan},~\bfnm{Dapeng}\binits{D.}}
(\byear{2008}).
\btitle{Reversibility of chordal {SLE}}.
\bjournal{Ann. Probab.}
\bvolume{36}
\bpages{1472--1494}.
\bid{doi={10.1214/07-AOP366}, issn={0091-1798}, mr={2435856}}
\end{barticle}
%
\endbibitem

\end{thebibliography}
\end{document}